\input amstex
\documentstyle{amsppt}
\magnification=1200
\hcorrection{0.2in}
\vcorrection{-0.4in}

\define\op#1{\operatorname{#1}}
\define\sball#1#2{{\underline B_{#2}^{#1}}}
\define\svolsp#1#2{{\volume(\partial \underline B_{#2}^{#1})}}
\NoRunningHeads
\NoBlackBoxes
\topmatter
\title Collapsed Manifolds With Ricci Bounded Covering Geometry
\endtitle
\abstract We study collapsed manifolds with Ricci bounded covering geometry i.e.,
Ricci curvature is bounded below and the Riemannian universal cover is
non-collapsed or consists of uniform Reifenberg points. Via Ricci flows'
techniques, we partially extend the nilpotent structural results of
Cheeger-Fukaya-Gromov, on collapsed manifolds with (sectional curvature)
local bounded covering geometry, to manifolds with (global) Ricci bounded
covering geometry.
\endabstract
\author Hongzhi Huang, Lingling Kong\footnote{Supported partially by NSFC Grant 11671070,11201058  and FRFCU Grant 2412017FZ002. The second author would like to thank Capital Normal University for
a warm hospitality during his visit},
Xiaochun Rong\footnote{Supported partially by a research fund from Capital Normal University.\hfill{$\,$}}
\&
Shicheng Xu\footnote{Supported partially by NSFC Grant 11401398 and by Youth Innovative Research Team of Capital Normal University.\hfill{$\,$}\\ 2010 Mathematics Subject Classification. Primary 53C21, 53C23, 53C24.\hfill{$\,$}}
\endauthor

\address School of Mathematics Science, Capital Normal University, Beijing 100048, P.R.C.
\endaddress
\email hyyqsaax\@163.com
\endemail
\address School of Mathematics and Statistics, Northeast Normal University, Changchu, JL 130024, P.R.China
\endaddress
\email KongLL111\@nenu.edu.cn
\endemail
\address Mathematics Department, Rutgers University
New Brunswick, NJ 08903 USA
\endaddress
\email rong\@math.rutgers.edu
\endemail
\address Mathematics Department, Capital Normal University, Beijing,
P.R.C.
\endaddress
\email  shichxu\@foxmail.com
\endemail

\date
\enddate
\endtopmatter
\document

\vskip4mm

\head 0. Introduction
\endhead

\vskip4mm

A complete manifold $M$ is called $\epsilon$-collapsed, if the volume of any unit ball on $M$ is less than $\epsilon$; one may normalize $M$ with a bound on curvature. Collapsed manifolds with bounded sectional curvature, $|\op{sec}|\le 1$,
has been extensively studied by Cheeger-Fukaya-Gromov ([CFG], [CG1, 2], [Fu1-3], [Gr]), and the basic discovery is a structure on $M$, called a nilpotent structure, consisting of compatible local nilpotent almost symmetric structures whose orbits point to all collapsed directions of the underlying metric i.e., any short geodesic loop at a point is locally homotopy
non-trivial and locally homotopic to a short loop in the orbit at the point. The
existence of such a structure
has found important
applications in Riemannian geometry (cf. [Ro1]).

As showed by examples, structures on collapsed manifolds with bounded Ricci
curvature are much more complicated; e.g., sectional curvature may blow up
at some points (\cite{An}, \cite{GW}, \cite{HSVZ}) or even everywhere (\cite{Li}). Hence, a realistic goal is to first restrict to certain interesting class. Under various additional conditions, via smoothing techniques one may show the existence of a nearby metric of bounded sectional curvature, so as to conclude, by applying the above results, nilpotent structures on collapsed manifolds with Ricci curvature bounded below or in absolute value (\cite{DWY}, \cite{PWY}, \cite{NZ}, etc).

In \cite{Ro2}, nilpotent structures were (directly) constructed on collapsed manifolds with Ricci bounded local covering geometry. A compact $n$-manifold $M$
is said to have Ricci bounded (resp. Ricci bounded below) local $(\rho,v)$-covering geometry, if
$|\op{Ric}_M|\le n-1$ (resp. $\op{Ric}_M\ge -(n-1)$), and local rewinding volume of
$B(x,\rho)$ satisfies
$$\op{vol}(B(\tilde x,\rho))\ge v>0, \quad \forall\, x\in M,\tag I $$
or local $(\delta,\rho)$-covering geometry, if for any $x\in M$,
$$\text{$\tilde x$ is a $(\delta,\rho)$-Reifenberg point},\tag II$$
where $\pi: (\widetilde{B(x,\rho)},\tilde x)\to (B(x,\rho),x)$ denotes the (incomplete)
Riemannian universal cover, $\pi(\tilde x)=x$, $B(x,\rho)$ denotes the metric
$\rho$-ball centered at $x$,
and $\tilde x$ is $(\delta,\rho)$-Reifenberg point, if for all $0<r\le \rho$, the Gromov-Hausdorff distance,
$$d_{GH}(B(\tilde x,r),\b B^0(r))<\delta \cdot r,$$
where $\b B^0(r)$ denotes an $r$-ball in $\Bbb R^n$ (cf. \cite{CC2}). Clearly, (II) implies
(I) and the converse does not hold. We point it out that for maximally collapsed
manifolds with Ricci bounded below i.e. whose diameters are small, local covering geometry is equivalent to (global) covering geometry.

Note that local bounded covering geometry (I) is a necessary condition
for collapsed manifolds with Ricci curvature
bounded below to have nilpotent structures (\cite{Ro2}). Indeed,
if $|\op{sec}_M|\le 1$, then there exist constants
$\rho(n), \delta(n)>0$ such that
$$\op{injrad}\left(B\left(\tilde x,{\frac {\rho(n)}2}\right)\right)\ge \epsilon(n)>0,\quad \forall\,\, x\in M,$$
where $\op{injrad}\,(U)$ denotes
the infimum of the injectivity radii on $U$. This crucial property was
pointed in \cite{CFG} based on the Gromov's theorem on almost flat manifold
\cite{Gr}, and a simple proof independent \cite{Gr} is given in \cite{Ro2}.

The local construction of nilpotent structures in \cite{Ro2} is independent of \cite{CFG}, and thus gives an alternative and relatively simple construction when restricting to collapsed manifolds with bounded sectional curvature. This seems
to be the only alternative local construction since \cite{Gr}.

Our goal is to show the existence of nilpotent structures on collapsed manifolds with Ricci bounded below and local covering geometry (I) or (II). We
point it out that without an upper bound on Ricci curvature, it seems to be difficult
in carrying out the local construction from \cite{Ro2}.

In the present paper, we will restrict ourself in the case that the Riemannian
universal cover $\tilde M$ of $M$ is not collapsed or all points on $\tilde M$ are $(\delta,\rho)$-Reifenberg points. Our main tools are Ricci flows (\cite{Ha}), the Perelman's pseudo-locality (\cite{Pe}), the structures on
Ricci limit spaces (\cite{CC1,2}) and equivariant GH convergence (cf. \cite{FY}); we will obtain a nearby
metric with bounded sectional curvature, to which we are able to
apply \cite{Gr} and \cite{CFG} to get nilpotent structures.

We now begin to state main results of this paper.

\proclaim{Theorem A}{\rm (Maximally collapsed manifolds with Ricci bounded blow local covering geometry)} Given $n, v>0$, there exist constants, $\epsilon(n,v), C(n)>0$, such that if
a complete $n$-manifold $M$ satisfies
$$\op{Ric}_M\ge -(n-1),\quad \op{vol}(B(\tilde p,1))\ge v>0,\quad \op{diam}(M )<\epsilon(n,v),$$
then $M$ is diffeomorphic to an infra-nilmanifold $N/\Gamma$, where $\tilde p$ is a point in
the Riemannian universal cover of $M$, $N$ is a simply connected nilpotent Lie group,
and $\Gamma$ is a discrete subgroup of
$N\rtimes \op{Aut}(N)$ such that $[\Gamma: \Gamma\cap N]\le C(n)$.
\endproclaim

Theorem A implies Gromov's theorem on almost flat manifolds i.e.,
$|\op{sec}_M|\le 1$ and $\op{diam}(M)<\epsilon(n)$; because which
satisfies local bounded covering geometry.

Theorem A does not hold when removing the non-collapsing condition
on $\tilde M$; any compact simply connected manifold of non-negative Ricci
curvature serves as a counterexample.

\remark{Remark \rm 0.1} Theorem A substantially improves previous results
on maximally collapsed manifolds with bounded Ricci or Ricci bounded below,
and various additional conditions to apply smoothing techniques (\cite{DWY}, \cite{PWY}, \cite{BW}); e.g.,
the conclusion of Theorem A was asserted in \cite{BW} under the condition
that $\op{Ric}_M\ge -(n-1)$ and $(\tilde M,\tilde p)$ is close to $(\Bbb R^n,0)$
in the Gromov-Hausdorff distance.
\endremark

A map between two metric spaces, $f: X\to Y$, is called an $\epsilon$-Gromov-Hausdorff approaximation, briefly an $\epsilon$-GHA,
if $f$ is an $\epsilon$-isometry and $f(X)$ is $\epsilon$-dense in $Y$. A sequence of Riemannian $n$-manifolds converges to
a metric space $X$ in the Gromov-Hausdorff topology, $M_i@>\op{GH}>>X$, if and only if there is
an $\epsilon_i$-GHA from $M_i$ to $X$, $\epsilon_i\to 0$.

\proclaim{Theorem B}{\rm (Nilpotent fibrations)} Let $M_i@>\op{GH}>>N$ be a sequence of complete $n$-manifolds such that
$$\op{Ric}_{M_i}\ge -(n-1),\qquad \op{vol}(B(\tilde p_i,1))\ge v>0,\quad \tilde p_i\in \tilde M_i,$$
where $N$ is a compact Riemannian manifold. Then for $i$ large, there is a smooth fibration,
$f_i: M_i\to N$, such that

\noindent (B1) An $f_i$-fiber is diffeomorphic to an infra-nilmanifold.

\noindent (B2) $f_i$ is an $\epsilon_i$-GHA, and a $(1-\Psi(\epsilon_i|n))$-H\"older map, $\epsilon_i\to 0$, where $\Psi(\epsilon_i|n)$ denotes a function in $\epsilon_i$ such that
$\Psi(\epsilon_i|n)\to 0$ as $\epsilon_i\to 0$, while $n$ is fixed.
\endproclaim

Theorem B can be viewed as a parametrized version of Theorem A; as
Theorem 1.2 a parametrized version of Theorem 1.1 (see Section 1).
Note that Theorem B will be false if one removes the non-collapsed
condition on $\tilde M_i$ (cf. [An]).

By (B2), the extrinsic diameter of any $f_i$-fiber is less than $\epsilon_i$,
which may not imply that the intrinsic diameter of an $f_i$-fiber is small
(that holds in case sectional curvature is bounded in absolute value, see
(1.2.4)). This is because in our approach to Theorem B, we use Ricci flows
to obtain a nearby collapsed metric with bounded sectional curvature,
by which we are able to apply the nilpotent structure result in \cite{CFG}.
Unfortunately, the Ricci flow metric is only weakly close
to the original metric i.e., their distance functions are bi-H\"older close.
Nevertheless, the following property holds (which is not used in the present paper): the normal subgroup, $\Lambda_i=\op{Im}[\pi_1(f_i\text{-fiber})\to \pi_1(M,p_i)]$,
preserves each component $\hat N_i$ of $\pi_i^{-1}(f_i\text{-fiber})$, in which
the $\Lambda_i$-orbit in $\hat N_i$ is $\epsilon_i$-dense.

The bundle projection map in Theorem B defines a pure nilpotent structure on $M$ (\cite{CFG}).
A pure nilpotent structure $\Cal N$ on an open connected subset $U$ of Riemannian manifold is defined by an $O(n)$-invariant
smooth affine fiber bundle on the orthogonal frame bundle over $U$, $\tilde f: (F(U),O(n))\to (Y,O(n))$, with fiber a nilpotent manifold (thus $Y$ is smooth) i.e., the $O(n)$-action preserves
both fibers and the structure group. The $O(n)$-invariance implies that $\tilde f$ descends to a map,
$f: U\to X=Y/O(n)$, such that the following diagram commutes:
$$\CD (F(U),O(n))@>\tilde f>>(Y,O(n))\\
@ V \pi  VV @ V \bar \pi VV\\
U@>f >> X=Y/O(n).
\endCD$$
Note
that a pure nilpotent structure, $f: U\to X$, is in general a singular fibration i.e. a singular fiber occurs when an $\tilde f$-fiber and $O(n)$-orbit meets on
a subset of positive dimension.

\proclaim{Theorem C} {\rm(Singular nilpotent fibrations)} Given $n$, there exists a small 
constant $\delta(n)>0$, such that if $M_i@>\op{GH}>>X$ is a sequence of complete $n$-manifolds converging to a compact space $X$ satisfying
$$\op{Ric}_{M_i}\ge -(n-1),\quad \text{$\tilde x_i$ is $(\delta,\rho)$-Reifenberg point},\quad 0<\delta\le \delta(n), \quad \forall\,\, \tilde x_i\in \tilde M_i,$$
then, for $i$ large, there is a singular fibration map, $f_i: M_i\to X$, such that

\noindent (C1) An $f_i$-fiber is diffeomorphic to an infra-nilmanifold.

\noindent (C2) $f_i$ is an $\Psi(\delta|n)$-GHA, and $(1-\Psi(\delta|n))$-H\"older map.
\endproclaim

A (mixed) nilpotent structure on a manifold $M$, $\Cal N=\{(U_i,\Cal N_i)\}$, consists of a locally finite open cover
for $M$, $\{U_i\}$, each $U_i$ admits a pure nilpotent structure $\Cal N_i$ such that the following compatible
condition holds: if $U_i\cap U_j\ne \emptyset$, then $U_i\cap U_j$ are both $\Cal N_i$-invariant and $\Cal N_j$-invariant such that $\Cal N_i$-orbits sit in an $\Cal N_j$-orbit or vice versa. By the compatibility, we define the $\Cal N$-orbit at a point by the $\Cal N_i$-orbit at the point of largest dimension.

\proclaim{Theorem D} \rm{(Mixed nilpotent structures)} Given $n,\rho>0$, there
exist constants, $\epsilon(n,\rho), \delta(n), C(n)>0$, such that if an $\epsilon$-collapsed compact
$n$-manifold, $\epsilon\le \epsilon(n,\rho)$, satisfies
$$\op{Ric}_M\ge -(n-1), \quad \text{$\tilde x$ is a $(\delta(n),\rho)$-Reifenberg point},\quad \forall\,\, \tilde x\in \tilde M,$$
then $M$ admits
a mixed nilpotent structure whose orbits at a point has extrinsic
diameter $<C(n)\epsilon$.
\endproclaim

\remark{Remark \rm 0.2} The existence of nilpotent structures on collapsed manifolds $M$ with $|\op{sec}_M|\le 1$ has found several important applications in Riemannian geometry, cf. \cite{Ro1}. By Theorems A-D, most of the applications, if not all, should hold on manifolds of Ricci bounded below and bounded covering geometry.
\endremark

We now briefly describe our approach to Theorems A-C: we will use Ricci flows (\cite{Ha}) and Perelman's pseudo-locality (\cite{Pe}, \cite{CTY}). Let
$g_i(t)$ denote the solution of
$$\frac {\partial g_i(t)}{\partial t}=-2\op{Ric}(g_i(t)),\quad g_i(0)=g_i.$$
If $g_i(t)$ exists for a definite time such that $|\op{Rm}(t^{-1}g_i(t))|\le 1$
and if $|d_{g_i}(x_{ik},x_{il})-d_{g_i(t)}(x_{ik},x_{il})|$ is much small than
$\sqrt t$ when either distance is less than $\sqrt t$ (see Lemma 1.11), then
we can apply nilpotent structural results in \cite{CFG} to obtain the desired fibration structures on $M_i$ for $i$ large.

We first claim that the assumptions in Theorems A and B guarantee that points
on the Riemannian universal cover $\tilde M_i$ are uniform Reifenberg points (see Lemma 2.1),
which in turn, implies that $\tilde M_i$ satisfies the isoperimetric inequality condition in Perelman's pseudo-locality theorem on Ricci flows (\cite{Pe}, \cite{CM}).
The verification of the claim is via equivariant convergence (\cite{FY}) and
structures on Ricci limit spaces (\cite{CC1,2}). Consequently, the Ricci flow solution $\tilde g_i(t)$ on $\tilde M_i$ exists for a definite time and
its curvature tensor, $|\op{Rm}(\tilde g_i(t))|\le \frac 1t$. Because the deck transformations
on $\tilde M_i$ are also isometries with respect to $\tilde g_i(t)$, $\tilde g_i(t)$ descends to
$g_i(t)$ on $M_i$ which is indeed the solution of Ricci flows on $M_i$.

In Theorem A, based on the local estimate on $|d_g-d_{g(t)}|$ (Lemma 1.11, cf. \cite{CRX})
we see that $\op{diam}(M,d_{t^{-1}g(t)})$ is small so that we may apply \cite{Gr}
to conclude Theorem A.

In the proof of Theorem B, if $d_{GH}((M_i,t^{-1}g_i(t)),t^{-\frac 12}N)$ is small, then we
may directly apply the fibration theorem in \cite {CFG} (cf. \cite{Fu1})
to conclude the desired result. Unfortunately, the estimate on local distance
functions is inadequate for the desired GH-closeness. Instead, we will show that
in our circumstances, the construction of fibration in \cite{Fu1} (cf. \cite{Ro1})
can be modified according to the local closeness of distance functions
(which does not require that $d_{GH}((M_i,t^{-1}g_i(t)),t^{-\frac 12}N)$ is small), see Theorem 2.2.

In the proof of Theorem C, we show that the limit space of the Ricci flow
metrics by a uniform definite time is bi-H\"older homeomorphic to $X$,
which is based on the bi-H\"older closeness of the initial metric and
its Ricci flow metric in \cite{BW} (cf. \cite{CRX}). Note that this also
implies a weak version of Theorem B.

To extend the approach of the present work to collapsed manifolds with Ricci bounded
blow and local $(\delta,\rho)$-bounded covering geometry, the key is to have
a type of Perelman's curvature estimate. We point it out that such an estimate,
if holds, would be essentially different from Perelman's pseudo-locality which
relies on local geometry. Precisely, around a point whose universal cover
satisfies the isopermetric inequality, there may not be a curvature bound
if there is a far away point whose local cover does not satisfy the
isoperimetric inequality (see below).

\example{Example 0.3}(Topping, cf. \cite{Cho}) Let $M_\epsilon$ be a topological sphere obtained by
capping an $\epsilon$-thin flat cylinder, $S^1(\epsilon)\times [-1,1]$, with two round $\epsilon$-hemispheres. With a slightly smoothing, we may assume that $M_\epsilon$
is rotationally symmetric Riemannian manifold of nonnegative sectional curvature.
For fixing $0<2\epsilon<\rho<<1$,
points that are $2\rho$-away from the centers of two hemispheres will satisfy local rewinding
Reifenberg points condition. However, the Ricci flow time is proportional to area of $M_\epsilon$ which is less than 2$\epsilon$. This may suggest that if one wants to have a type of Perelman's pseuo-locality estimate on curvature of flowed metric, one has to use information on
local rewinding volume at points far away.
\endexample

The rest of the paper is organized as follows:

In Section 1, we will supply notions and basic properties that will be used in the proofs of Theorems A-D.

In Section 2, we will prove Theorems A-D.

In Appendix, for convenience of readers, we will outline a proof of Theorem 1.2 via
embedding method in \cite{Ro1} with minor modifications, and we will present a proof for the
bi-h\"older estimate of distance functions of the original metric and the Ricci flows
with pseudo-locality (\cite{BW}).

\vskip4mm

\head 1. Preliminaries
\endhead

\vskip4mm

In this section, we will supply notions and basic results that will be
used through the rest of the paper.

\vskip4mm

\subhead a. N-structures and collapsing with  bounded sectional curvature
\endsubhead

\vskip4mm

We will briefly recall main structural results on collapsed manifolds
with bounded sectional curvature (\cite{CFG}, \cite{CG1,2}, \cite{Fu1-3}, \cite{Gr}).

Let $U$ be an open subset of an $n$-manifold $M$ and let
$P: F(U)\to U$ denote the orthogonal frame bundle on $U$. A pure N-structure, $\Cal N$, on
$U$ refers to an $O(n)$-invariant fiber bundle, $\tilde f: F(U)\to Y$,
such that an $\tilde f$-fiber is a nilmanifold, $N/\Gamma$, with $N$ a
simply connected nilpotent Lie group, $\Gamma$ a co-compact lattice of
$N$, and the structural group is a subgroup of $N\rtimes \op{Aut}(N)$.
The $O(n)$-invariance implies that the $O(n)$-action on $F(U)$ descends
to an $O(n)$-action on $Y$, and $\tilde f$ descends to a map,
$f: U\to Y/O(n)$, such that the following diagram commutes,
$$\CD (F(U),O(n))@>\tilde f>> (Y,O(n))\\
@V P VV   @ V P VV \\
U@>f>> Y/O(n).
\endCD$$
Note that each $f$-fiber is an infra-nilmanifold which may not have
constant dimension; when $O(n)$-orbits and $\tilde f$-fiber intersect more
than isolated points. The minimal dimension of all $f$-fibers is called
the rank of the pure N-structure $\Cal N$. If $U=M$, we say that $M$
admits a pure N-structure. A subset of $U$ is called invariant if $U$ is
the disjoint union of $f$-fibers (or $f$-orbits), and a metric on $U$ is
called invariant, if the induced metric on every $f$-fiber is left invariant.

A mixed N-structure on $M$ is defined by a locally finite open cover,
$\{U_i\}$, for $M$, together with a pure N-structure $\Cal N_i$ on $U_i$,
which satisfies the following compatibility condition: whenever
$U_i\cap U_j\ne \emptyset$, $U_i\cap U_j$ are both invariant with respect
to $\Cal N_i$ and $\Cal N_j$, and on $F(U_i\cap U_j)$, every $\tilde f_i$-fiber
is a union of $\tilde f_j$-fibers or vice versa.

An important case of a mixed N-structure is that all $\tilde f_i$-fibers are tori, which
is called an F-structure on $M$. In [CG1], Cheeger-Gromov
showed that if a complete manifold admits a (mixed) F-structure of positive rank,
then $M$ admits a  one-parameter family of invariant metrics, $g_\epsilon$,
which are $\epsilon$-collapsed metrics such that $|\op{sec}_{g_\epsilon}|\le 1$
and every orbit collapse to a point.  Cheeger-Gromov proposed a similar
construction with respect to a given mixed N-structure  of positive rank,
this was carried out in \cite{CR}.

We now begin to recall main structural results on collapsed manifolds with bounded sectional curvature.

\proclaim{Theorem 1.1}  {\rm(\cite{Gr}, \cite{Ru})} There exists a constant $\epsilon(n)>0$
such that if a compact $n$-manifold $M$ satisfies
$$|\op{sec}_M|\le 1,\quad \op{diam}(M)<\epsilon(n),$$
then $M$ is diffeomorphic to an infra-nilmanifold.
\endproclaim

The converse of Theorem 1.1 holds (\cite{Gr}).  The following can be viewed as a bundle
version of maximal collapsed manifolds.

\proclaim{Theorem 1.2} {\rm(\cite{Fu1}, \cite{CFG})} Given $n$, there exist constants,
$\epsilon(n), c(n)>0$, such that if a compact $m$-manifold $N$ and an $n$-manifold $M$ satisfy $$|\op{sec}_M|\le 1, \quad |\op{sec}_N|\le 1,\quad \op{injrad}(N)\ge 1,\quad  d_{GH}(M,N)<\epsilon\le \epsilon(n),$$
then there is a fibration map, $f: M\to N$, such that

\noindent (1.2.1) $f$ is a $C\cdot \epsilon$-GHA, where $C$ is a constant.

\noindent (1.2.2) An $f$-fiber is diffeomorphic to an infra-nilmanifold.

\noindent (1.2.3) $f$ is $\epsilon$-Riemannian submersion, $e^{-\Psi(\epsilon|n)}\le \frac{|df_x(\xi)|}{|\xi|}\le e^{\Psi(\epsilon|n)}$, where $\xi$ is orthogonal to the $f$-fiber at $x$.

\noindent (1.2.4) The second fundamental form of $f$-fibers,
$|II_f|\le c(n)$.
\endproclaim


\remark{Remark \rm 1.3} In \cite{CFG}, $\Psi(\epsilon|n)=c(n)\sqrt \epsilon$, and
(1.2.1) may not been seen; due to the fact that averaging operation is performed
on a $\sqrt \epsilon$-ball. In \cite{Fu}, \cite{Ya} (cf. \cite{Ro1}), $f$ is
constructed via embedding $N$ and $M$ into an Euclidean space via (averaging)
distance functions of two $d_{GH}$-close $\epsilon$-nets on $N$ and $M$. The
embedding method works with condition, $\op{sec}_M\ge -1$, where $f$ is
$C^1$-smooth (so (1.2.4) may not be seen). We point it out that
with a minor modification of the construction of $f$ in \cite{Ro1},
plus the (additional) condition $\op{sec}_M\le 1$, one gets that
$f$ is smooth and (1.2.4). For convenience of readers, we give a brief
proof in Appendix.
\endremark

For a collapsing sequence, $M_i@>\op{GH}>>X$, with $|\op{sec}_{M_i}|\le 1$ and
$X$ a compact metric space, the associate sequence of frame bundles equipped with
canonical metrics, $F(M_i)@>\op{GH}>>Y$, and $Y$ is a manifold ([Fu2]). By extending
Theorem 1.2 to an $O(n)$-equivariant version on $F(M_i)$, one gets a singular fibration on
$M_i$.

\proclaim{Theorem 1.4} {\rm(\cite{Fu2}, \cite{CFG})} Let $M_i@>\op{GH}>>X$ be a sequence
of compact $n$-manifolds with $X$ a compact length space and $|\op{sec}_{M_i}|
\le 1$.  Then for $i$ large, there is a singular fibration map, $f_i: M_i\to X$,
such that

\noindent (1.4.1) An $f_i$-fiber is diffeomorphic to an infra-nilmanifold.

\noindent (1.4.2) $f_i$ is an $\epsilon_i$-GHA, $\epsilon_i\to 0$.

\noindent (1.4.3) $f_i$ is a $C(n)\epsilon_i$-submetry i.e., for all $x_i\in M_i$ and $r>0$,
$$B(f(x_i),e^{-C(n)\epsilon_i}r)\subseteq f_i(B(x_i,r))\subseteq B(f_i(x_i),e^{C(n)\epsilon_i} r).$$
\endproclaim

Note that the above construction of $f_i$ can be made local, and thus a bound
on diameter can be removed for the singular fibration structure.

\proclaim{Theorem 1.5} {\rm(\cite{CG1,2}, \cite{CFG})} There exists a constant
$\epsilon(n)$ such that if a complete $n$-manifold $M$ satisfies
$$|\op{sec}_M|\le 1,\quad \op{injrad}(M,x)<\epsilon(n),\quad \forall\, x\in M,$$
then $M$ admits a (mixed) N-structure of positive rank and at any $x\in M$, the
orbit points to all collapsed directions with respect to a fixed scale.
\endproclaim

\vskip4mm

\subhead b. Ricci flows and the Perelman's pseudo-locality
\endsubhead

\vskip4mm

Given a complete manifold $(M,g)$, the Ricci flow is the solution of
$$\frac {\partial g(t)}{\partial t}=-2\op{Ric}(g(t)),\quad g(0)=g.$$
If $\sup |\op{sec}_g|\le 1$ (e.g., $M$ is compact), then there is a short time complete smooth solution (\cite{Ha}, \cite{Sh1,2}) which, for $t\in (0,T(n)]$, satisfies $$\sup_M|\nabla^m\op{Rm}(g(t))|^2\le \frac {C_m}{t^m},\quad  m=0, 1, ..., $$
By [CZ], smooth solutions with bounded curvature tensor is unique. Let $T_{\max}$ denote
the supremum of $t$ for $g(t)$ exists. If $T_{\max}<+\infty$, then $\sup |\op{Rm}(g(t))|\to \infty$
as $t\to T_{\max}$. If a group $G$ acts isometrically with respect to $g(0)$, then
$G$ acts by isometries with respect to $g(t)$.

In Riemannian geometry, Ricci flow has been a powerful tool to raise the regularity
of original metrics, which requires a definite flow time. A fundamental result
for estimating a lower bound on $T_{\max}$ is the following Perelman's
pseudo-locality theorem on Ricci flows (\cite{Pe}).

\proclaim{Theorem 1.6}{\rm (Perelman's pseudo-locality)}
Given $n\ge2,\alpha>0$, there exist $\epsilon_0$ and $\delta$ satisfying the following. Let $(M^n,g(t))$ be a complete solution of the Ricci flow with bounded curvature, where $t\in[0,(\epsilon r_0)^2],$ $\epsilon\le\epsilon_0$, $r_0\in (0,\infty)$, and let $x_0\in M$ be a point such that scalar curvature $$\op{R}_{g(0)}(x)\ge-r_0^{-2},\quad \text{for all }x\in B_{g(0)}(x_0,r_0),$$and $B_{g(0)}(x_0,r_0)$ is $\delta$-almost isoperimetrically Euclidean: $$\left(\op{vol}_{g(0)}(\partial\Omega)\right)^{n}\ge (1-\delta)c_n\left(\op{vol}_{g(0)}(\Omega)\right)^{n-1}$$ for every regular domain $\Omega\subset B_{g(0)}(x_0,r_0),$ where $c_n$ is the Euclidean isoperimetric constant. Then
$$|\op{Rm}|(x,t)\leq\frac{\alpha}{t}
+\frac{1}{(\epsilon_0r_0)^2},\quad\op{vol}_{g(t)}(B_{g(t)}(x,\sqrt{t}))\ge C(n)\sqrt{t}^n$$ for all $x\in M$ such that $d_{g(t)}(x,x_0)<\epsilon_0r_0$ and $t\in(0,(\epsilon r_0)^2]$, where $C(n)$ is a constant depending on $n$.
\endproclaim

\remark{Remark \rm 1.7} Theorem 1.6 implies that the complete solution of the Ricci flow exists at least for $t\in[0,(\epsilon_0r_0)^2]$, provided that the curvature conditions and the
isoperimetric inequality hold everywhere.
\endremark

The original statement of Theorem 1.6 in \cite{Pe} restricts to compact $M$.
For a complete noncompact $M$, see \cite{CTY} (cf. Chapter 21 in \cite{Cho}).

In case that $\op{Ric}\ge-(n-1)$, the almost isoperimetrically Euclidean condition is equivalent to the uniformly Reifenberg condition. For our purpose, we will formulate a special case of Corollary 1.3
in \cite{CM} as follows:

\proclaim{Theorem 1.8}{\rm (\cite{CM})}
Given $n\ge2$, there exist $\epsilon(n),\delta(n),C(n)>0$ such that, if $M$ is an $n$-manifold with $\op{Ric}_M\ge-(n-1)\delta$ and $d_{GH}(B(p,1),\b{\rm{B}}^0(1))\le\delta\le\delta(n)$ where $B(p,1)\subset M$ is relative compact, then given any $\epsilon\le\epsilon(n)$, for every regular domain $\Omega\subset B(p,\epsilon)$, the following almost Euclidean isoperimetric inequality holds, $$\op{vol}(\partial\Omega)^n\ge(1-C(n)(\epsilon+\delta))c_n\op{vol}(\Omega)^{n-1},$$ where $c_n$ is given in Theorem 1.6.
\endproclaim

We will apply Theorem 1.6 and 1.8 in the following situation: a compact
$n$-manifold $M$ satisfies that
$$\op{Ric}_M\ge -(n-1), \quad d_{GH}(B(x,1),\b B^0(1))<\delta.$$
Scaling $M$ by $\delta^{-\frac 12}>1$, $\delta^{\frac 12}\le \delta(n)$ in Theorem 1.8, we obtain
$$\op{Ric}_{\delta^{-\frac 12}M}\ge -(n-1)\delta,\quad d_{GH}(B_{\delta^{-\frac 12}M}(x,\delta^{-\frac 12}),\b B^0(\delta^{-\frac 12}))<\delta^{\frac 12}.\tag 1.9$$
By Theorem 1.8, we conclude that any $\Omega$ in a unit ball
on $\delta^{-\frac 12}M$ satisfies
$$\op{vol}(\partial\Omega)^n\ge(1-C(n)(\delta^{\frac 12}+\delta))c_n\op
{vol}(\Omega)^{n-1}.$$
Because the above inequality is scaling invariant, $M$ also satisfies the above isoperimetric
inequality condition. In short, we are able to apply Theorem 1.6 under condition
(1.9).

An important consequence of the pseudo-locality is a local distance estimate (see below). Let the assumption be as in Theorem 1.6 which holds everywhere (see Remark 1.7), and for the sake of simple notation, let $d_t
:=d_{g(t)}$. Then by Hamilton's integral version of Myer's theorem (cf. Lemma 8.3 in \cite{Pe}), for all $x,y\in M$ and $0\le t\le (\epsilon r_0)^2$,
 $$d_{0}(x,y)\le d_{t}(x,y)+c(n)\sqrt{\alpha t}.\tag{1.10}$$
Conversely, as observed in \cite{BW} (cf. see Lemma 2.10 \cite{CRX}),
$d_t(x,y)$ is also controlled by $d_0$ in the scale of $\sqrt{t}$. For our
purpose, we need the following modification.

\proclaim{Lemma 1.11} Let the assumptions be as in Theorem 1.6.
For any $x,y\in M$ with $d_t(x,y)\leq \sqrt{t}$ or
$d_0(x,y)\le \sqrt t$, we have
$$|d_t(x,y)-d_0(x,y)|\leq \Psi(\alpha| n) \sqrt{t}, \tag 1.12$$
where $\Psi(\alpha|n)\to 0$ as $\alpha\to 0$.
\endproclaim

\demo{Proof}
By (1.10), it suffices to show $d_0(x,y)\ge d_t(x,y)-\Psi(\alpha|n)\sqrt{t}$. This holds
if $d_t(x,y)\leq \sqrt{t}$, see Lemma 2.10 in \cite{CRX} where
$d_t(x,y)<\sqrt t$ after normalizing $\rho(n)=1$ ($\rho(n)<1$), but it actually
holds without a normalization (see below).

Assume that $d_{0}(x,y)\le \frac 12c(n)\sqrt{t}$, and we will show that $d_{0}(x,y)\le \frac 12c(n)\sqrt{t}$
implies that $d_{t}(x,y)\le c(n)\sqrt{t}$.

Arguing by contradiction, let $\gamma$ be a segment from $x$ to $y$ at time $0$ such that $d_t(x,y)>c(n)\sqrt{t}$. Then in time $t$ we can find a point $z\neq y$ on $\gamma$ such that $d_{t}(x,z)=c(n)\sqrt{t}$. By (1.12) for the case $d_t(x,y)\le c(n)\sqrt{t}$,
$$d_0(x,z)\ge d_t(x,z)-\Psi(\alpha|n)\sqrt{t}.$$ So
for sufficient small $\alpha$ such that $\Psi(\alpha|n)<\frac{1}{2}c(n)$,
$$d_{0}(x,y)>d_{0}(x,z)\ge \frac{1}{2}c(n)\sqrt{t},$$
a contradiction.

Finally, we point that $|d_t(x,y)-d_0(x,y)|<\Psi(\alpha|n)\sqrt t$ with
$d_t(x,y)<c(n)\sqrt t$ or $d_0(x,y)<c(n)\sqrt t$ implies that it holds with $d_t(x,y)<\sqrt t$ or $d_0(x,y)<\sqrt t$; dividing $[0,1]$ by $c(n)<1$ and on each subinterval applying that $|d_t(x,y)-d_0(x,y)|<\Psi(\alpha|n)\sqrt t$.
\qed\enddemo

\remark{Remark \rm 1.13} Let $(M,g)$ be a compact $n$-manifoild, and let $g(t)$ denote the 
unique short time Ricci flow solution. Then the pullback $\tilde g(t)$ on $\tilde M$ is 
also a Ricci flow solution of bounded curvature. If $(\tilde M,\tilde g)$ satisfies that $\op{Ric}_{\tilde g}\ge -(n-1)$ and all points are $(\delta(\alpha|n),\rho)$-Reifenberg, then Theorems 1.6 applies to $\tilde g(t)$, and thus $g(t)$ satisfies the same curvature estimate. Therefore,  Ricci flow solution on $(M,g)$ exists in $[0, T(n,\rho,\alpha)]$, and
$|\op{Rm}(g(t))|\le \frac{\alpha}{t}, t\in (0,T(n,\rho,\alpha)]$.

Moreover, the local distance estimate (1.12) also holds on $(M,g)$.
\endremark

\vskip4mm

\subhead c. Ricci limit spaces
\endsubhead

\vskip4mm

A pointed metric space $(X,p)$ is called a Ricci limit space, if there is a sequence of complete $n$-manifolds of $\op{Ric}_{M_i}\ge -(n-1)$ such that $(M_i,p_i)@>\operatorname{GH}>>(X,p)$. A Ricci limit space $(X,p)$ is called $v$-non-collapsed, if $\op{vol}(B(p_i,1))\ge v>0$. Let $\Cal M(n,-1)$ (resp.
$\Cal M(n,-1,v)$) denote the collection of Ricci limit space (resp. $v$-non-collapsed Ricci limit spaces). We will briefly recall some basic results in the Cheeger-Colding theory on Ricci limit spaces that will be used in our proof of Theorem A-D.

For any $q\in X\in \Cal M(n,-1)$, and any $r_i\to \infty$, by Gromov's compactness it is easy to see, passing to a subsequence, $(r_iX,q)@>\op{GH}>>(C_qX,o)$, which is called a tangent cone of $X$ at $q$. A tangent cone
may depend on the choice of a subsequence. A point $q$ is called regular, if its tangent cone
is unique and is isometric to an Euclidean space. A point is called singular, if it is not
regular.

\proclaim{Theorem 1.14} {\rm(\cite{Co1}, \cite{CC1, 2}) Let $(X,p)\in \Cal M(n,-1)$.

\noindent (1.14.1) If a tangent cone, $C_qX$, contains a line, then $C_qX$ splits off an $\Bbb R$-factor.

\noindent (1.14.2) The set of regular points in $X$ is dense.

\noindent (1.14.3) If $X\in \Cal M(n,-1,v)$, then $C_qX$ is an $n$-dimensional metric cone
i.e., $C_qX\cong \Bbb R^k\times C(Z)$, where $C(Z)$ is a metric cone over a length space $Z$ of $\op{diam}(Z)<\pi$.
In particular, any regular point has a tangent cone $\Bbb R^n$.
\endproclaim

\proclaim{Theorem 1.15}{\rm([Co2])} For each $\epsilon,n,v$, there exists $\delta(\epsilon,n)$ and $r_0(\epsilon,n)$ satisfying the follows. Given $(X,p)\in \Cal M(n,-1,v)$ and $q\in X$, if $d_{GH}(B(q,r_0),\b{\rm{B}}^0(r_0))\le r_0\delta$, then $q$ is an $(\epsilon,r_0)$-Reifenberg point.
\endproclaim

By Colding's volume convergence ([Co2]) and Theorem 1.15, for any $x\in B(q,r_0)$ with distance $d(q,x)<(1-\epsilon)r_0$, $x$ is $(\epsilon,s)$-Reifenberg, where $s=(1-\epsilon)r_0-d(q,x)$ (cf. [CC2]).
\vskip4mm

\subhead d. Equivariant Gromov-Hausdorff convergence
\endsubhead

\vskip4mm

The reference of this part is [FY] (cf. [Ro1]).

Let $X_i@>\op{GH}>>X$ be a convergent sequence of compact length metric spaces, i.e.,
there are a sequence $\epsilon_i\to 0$ and a sequence of $\epsilon_i$-GHA maps,
$h_i: X_i\to X.$
Assume that $X_i$ admits a closed group $\Gamma_i$-action by isometries. Then
$(X_i,\Gamma_i)@>\op{GH}>>(X,\Gamma)$ means that there are a sequence $\epsilon_i\to 0$ and
a sequence of $(h_i,\phi_i,\psi_i)$, $h_i: X_i\to X$, $\phi_i:
\Gamma_i\to \Gamma$ and $\psi_i: \Gamma\to \Gamma_i,$ which are $\epsilon_i$-GHAs
such that for all $x_i\in X_i, \gamma_i\in \Gamma_i$ and $\gamma\in \Gamma$,
$$d_X(h_i(x_i), \phi_i(\gamma_i)h_i(\gamma_i^{-1}(x_i)))<\epsilon_i,\quad
d_X(h_i(x_i),\gamma^{-1}(h_i(\psi_i(\gamma)(x_i)))<\epsilon_i,\tag 1.16$$
where $\Gamma$ is a closed group of isometries on $X$, $\Gamma_i$ and
$\Gamma$ are equipped with the induced metrics from $X_i$ and $X$.
We call $(h_i,\phi_i,\psi_i)$
an $\epsilon_i$-equivariant GHA.

When $X$ is not compact, then the above notion of equivariant convergence
naturally extends to a pointed version $(h_i,\phi_i,\psi_i)$: $h_i:
B(p_i,\epsilon_i^{-1})\to B(p,\epsilon_i^{-1}+\epsilon_i)$, $h_i(p_i)=p$,
$\phi_i: \Gamma_i(\epsilon_i^{-1})\to \Gamma(\epsilon_i^{-1}+\epsilon_i)$,
$\phi_i(e_i)=e$, $\psi_i: \Gamma(\epsilon_i^{-1})\to \Gamma_i(\epsilon_i^{-1}
+\epsilon_i)$, $\psi_i(e)=e_i$, and (1.16) holds whenever the multiplications
stay in the domain of $h_i$, where $\Gamma_i(R)=\{\gamma_i\in\Gamma_i,\,
\,d_{X_i}(p_i,\gamma_i(p_i))\le R\}$.

\proclaim{Lemma 1.17}  Let $(X_i,p_i)@>\op{GH}>>(X,p)$, where $X_i$ is a  complete and
locally compact length space. Assume that $\Gamma_i$ is a closed group of isometries
on $X_i$. Then there is a closed group $G$ of isometries on $X$
such that passing to a subsequence, $(X_i,p_i,\Gamma_i)@>\op{GH}>>(X,p,G)$.
\endproclaim

\proclaim{Lemma 1.18}  Let $(X_i, p_i, \Gamma_i)@>\op{GH}>>(X,p,G)$, where $X_i$ is
a complete, locally compact length space and $\Gamma_i$ is a closed subgroup
of isometries. Then $(X_i/\Gamma_i,\bar p_i)@>\op{GH}>>(X/G,\bar p)$.
\endproclaim

For $p_i\in X_i$, let $\Gamma_i=\pi_1(X_i,p_i)$ be the fundamental group.
Assume that the universal covering space, $\pi_i: (\tilde X_i,
\tilde p_i)\to (X_i,p_i)$,
exists.
By Lemmas 1.17 and 1.18, we have the following commutative diagram:
$$\CD (\tilde X_i,\tilde p_i,\Gamma_i)@>\op{GH}>>(\tilde X,\tilde p,G)
\\ @VV \pi_i V   @VV \pi  V \\
(X_i,p_i)@>\op{GH}>> (X,p)=(\tilde X/G,p).\endCD \tag 1.19$$

\vskip4mm

\head 2. Proof of Theorems A-D
\endhead

\vskip4mm

\subhead a. Proof of Theorem A and B
\endsubhead

\vskip4mm
\demo{Proof of Theorem A}

By the Gromov's compactness, it suffices to show that given a sequence of compact $n$-manifolds $M_i$
satisfying
$$M_i@>\operatorname{GH}>>\{p\},\quad \op{Ric}_{M_i}\ge -(n-1),\quad \op{vol}(B(\tilde p_i,1))\ge v>0,$$
then for $i$ large, $M_i$ is diffeomorphic to an infra-nilmanifold. We will
show that $M_i$ admits a metric satisfying the conditions of Theorem 1.1.

Consider the equivariant convergence (see (1.19)), $$\CD (\tilde M_i,\tilde p_i,\Gamma_i)@>\operatorname{GH}>>(\tilde X,\tilde p,G)\\ @ V \pi_i VV @ V \pi VV\\(M_i,p_i)@>\operatorname{GH}>>\{p\}=\tilde X/G. \endCD$$
Since $\op{vol}(B(\tilde p_i,1))\ge v>0$, $\dim (\tilde X)=n$.
Because the limit group $G$ acts transitively on $\tilde X$ which contains a regular point
(see (1.14.2)), all points in $\tilde X$ are regular.
By Theorem 1.15, it is easy to see that given $\epsilon>0$, up to a
rescaling by constant, all points on $\tilde M_i$ are $(\epsilon,1)$-Reifenberg points, for $i$ large. To apply Theorem 1.6,
we specify $\alpha<1$ such that $\Psi(\alpha|n)<\frac {\epsilon(n)}2$, where
$\epsilon(n)$ is given by Theorem 1.1 and $\Psi(\alpha|n)$ is given in Lemma 1.11.
By Theorem 1.6, a complete solution $\tilde g_i(t)$ of Ricci flow on $\tilde M_i$,
$$\frac{\partial \tilde g_i(t)}{\partial t}=-2\operatorname{Ric}(\tilde g_i(t)), \quad \tilde g_i(0)=\tilde g_i,$$
exists on $[0,T(n,\alpha)]$, where the initial metric $\tilde g_i$ is
the pullback metric on $\tilde M_i$ from $M_i$. Because $\Gamma_i$ acts
on $\tilde M_i$ by isometries with respect to $\tilde g_i(t)$, we use $g_i(t)$ to
denote the quotient metric of $\tilde g_i(t)$ on $M_i$, which is
indeed the Ricci flow solution on $(M_i,g_i)$.
By Lemma 1.11, we derive that for any fixed $t_0\in (0,T(n,\alpha)]$ and  $i$ large,
$$\split d_{GH}((M_i,&t_0^{-1}g_i(t_0)),\{p\})\\
\leq~ & d_{GH}((M_i,t_0^{-1}g_i(t_0)),(M_i,t_0^{-1}g_i))+d_{GH}((M_i,t_0^{-1}g_i),\{p\})
\\ \leq~&\Psi(\alpha|n)+\epsilon_it_0^{-1/2}<\epsilon(n),\endsplit$$
where $\epsilon_i=d_{GH}((M_i,g_i),\{p\})$ i.e., $t_0^{-1}g_i(t_0)$ is a metric on
$M_i$ satisfying the conditions of Theorem 1.1.
\qed\enddemo

Our approach to Theorem B is also by smoothing technique via Ricci flows. As seen
in the proof of Theorem A, the first step is to check that all points on
the equivariant limit space are regular.

\proclaim{Lemma 2.1} Let $M_i @>\operatorname{GH}>> N$ be a sequence of compact $n$-manifolds satisfying Theorem B. Consider the following commutative diagram,
$$\CD (\tilde{M}_i,\tilde{p}_i,\Gamma_i)@>\op{GH}>>(Y,\tilde p,\Gamma)\\
@ V \pi_i  VV @ V \pi VV\\
(M_i,p_i)@>\operatorname{GH}>> (N,p)=Y/\Gamma,
\endCD$$
Then any point in $Y$ is regular.
\endproclaim
\demo{Proof}
Let us assume a sequence of $s_i\to \infty$ such that  $(s_i N, p)@>\operatorname{GH}>> (\Bbb R^m,0)$ and $(s_i Y, \tilde p, \Gamma) @>\operatorname{GH}>> (C_{\tilde p}(Y),\tilde p_*,  G),$ by Theorem 1.14,
we get the following commutative diagram,
$$\CD (s_iY,\tilde p,\Gamma)@>\operatorname{GH}>>(\Bbb R^k\times C(Z),(0,y_\infty),G)\\
@ V \pi_i  VV @ V \pi VV\\
(s_iN,p)@>\operatorname{GH}>> (\Bbb R^m,0),
\endCD$$
where $(C(Z),y_\infty)$ is a metric cone and contains no line. It suffices to show that $C(Z)$ is a point, which is equivalent to say that for any point $(x,y)\in \Bbb R^k \times C(Z)$, $y$ must coincide with $y_\infty$.

Argue by contradiction, let us assume that there is another point $y(\neq y_\infty) \in C(Z)$. Then $C(Z)$ contains at least one ray but contains no line. So $G$ preserves the subset $\Bbb R^k\times \{y_\infty\}$,  and thus $\pi(x,y)\neq \pi(x',y_\infty)$ for any $x,x'\in \Bbb R^k$.  In particular, $\pi(x,y)\neq \pi(0,y_\infty)$. Let $\gamma$ be a line in
$\Bbb R^m$ connecting $\pi(0,y_\infty)$ and $\pi(x,y)$, and let
$\tilde{\gamma}$ be the lifted line in $\Bbb R^k\times C(Z)$ passing through $(0,y_\infty)$ such that $\pi(\tilde \gamma)=\gamma$. Because $C(Z)$ contains no line,  $\tilde\gamma(t)$ must have the form $(\omega(t),y_\infty)$, which implies that $\tilde \gamma$ is a line in $\Bbb R^k\times \{y_\infty\}$. Then for any $t$,
$\gamma(t)=\pi(\tilde{\gamma}(t))=\pi(\omega(t),y_\infty)\neq \pi(x,y)$. On the other hand, by the construction of $\gamma$, there is some $t_0$ such that $\gamma(t_0)=\pi(x,y)$, a contradiction.
\qed\enddemo

Let $(M_i,g_i)@>GH>>(N,g_N)$ be as in Theorem B. By Lemma 2.1 and
the discussion following Theorem 1.8, for $i$ large we may apply Theorem 1.6 to conclude
that for any $0<\alpha\le 1$, the Ricci flow on $(M_i,g_i)$, $g_i(t)$,
exists for $t\in(0, T(n,\alpha)]$ such that $|\op{Rm}(t^{-1}g_i(t))|\le \alpha$
and for $x_i,y_i\in M_i$ with $d_{t^{-1}g_i(t)}(x_i,y_i)\le 1$,
$$|d_{t^{-1}g_i(t)}(x_i,y_i)-d_{t^{-1}g_i}(x_i,y_i)|\leq \Psi(\alpha|n).$$

Observe that if $d_{GH}((M_i,T^{-1}g_i(T)),(N,T^{-1}g_N))<\epsilon(N)$ with $T=T(n,\alpha)$, then Theorem B follows from Theorem 1.2.

Because the above global GH-closeness may not be true, we will need the following
generalization of Theorem 1.2, which says that Theorem 1.2 holds, if there is
an $(1,\epsilon)$-GHA (here `$1$' is a normalization) i.e., a map, $\phi: M\to N$,
which is $\epsilon_i$-onto and restricting to every unit ball, $\phi$ is an
$\epsilon$-isometry, and $\op{diam}(\phi^{-1}(x))<\epsilon$. Note that an $\epsilon$-GHA is an $(1,\epsilon)$-GHA,
but the converse may not hold.

\proclaim{Theorem 2.2} Given $n$, there exist constants
$\epsilon(n), c(n)>0$, such that if $\phi: M\to N$ is an $(1,\epsilon)$-GHA from a
compact $n$-manifold $M$ to a compact $m$-manifold $N$ satisfying
$$|\op{sec}_M|\le 1,\quad |\op{sec}_N|\le 1,\quad \op{injrad}(N)\ge 1,\quad 0<\epsilon\le \epsilon(n),$$
then there is a fibration map, $f: M\to N$, such that

\noindent (2.2.1) $f$ is an $(\frac 12,C\epsilon)$-GHA, where $C$ is a constant.

\noindent (2.2.2) An $f$-fiber is diffeomorphic to an infra-nilmanifold.

\noindent (2.2.3) $f$ is an $\epsilon$-Riemannian submersion, $e^{-\Psi(\epsilon|n)}\le \frac{|df_x(\xi)|}{|\xi|}\le e^{\Psi(\epsilon|n)}$, where $\xi$ is orthogonal to the $f$-fiber at $x$.

\noindent (2.2.4) The second fundamental form of $f$-fibers,
$|II_f|\le c(n)$.
\endproclaim

Note that the construction of nilpotent fibration in \cite{CFG} is local
and is completely determined by local geometry of a $C^1$-close metric
of higher regularity. Clearly, the argument in \cite{CFG} goes through with $(1,\epsilon(n))$-GHA;
if one only assumes that $\phi: M\to N$ is $\epsilon$-GHA on any unit ball and $\epsilon$-onto, then the argument in \cite{CFG} works through without any
modification, except that a fiber may not be connected; so the connectedness
of fibers is guaranteed by that $\op{diam}(\phi^{-1}(x))<\epsilon$.

Because Theorem 2.2 generalizes Theorem 1.2 (see Remark 1.3), we will also
briefly explain how Theorem 2.2 follows from the proof of Theorem 1.2 in
\cite{Ro1} (see Appendix), with an obvious modification in choosing a finite
$\epsilon$-dense subset $\{x_i\}\subset M$ (originally based on $d_{GH}(M,N)<\epsilon)$:
Given an $\epsilon$-net, $\{y_i\}$ of $N$, a consequence of an $(1,\epsilon)$-GHA $\phi: M\to N$ implies that  $\phi^{-1}(y_i)$ consists of `$\epsilon$-connected' components, each component
has diameter $<\epsilon$
and the distance of any two components $>>\epsilon$. We now choose, given each $y_i$, $x_{ij}
\in M$ such that $x_{ij}$ is a representative for each component.
Because the embedding construction depends only on local information,
the construction in \cite{Ro1} (cf. Appendix) goes through, with obvious
modifications corresponding to the additional condition, $\op{sec}_M\le 1$.

\demo{Proof of Theorem B}

Let $(M_i,g_i)@>GH>>(N,g_N)$ be as in Theorem B. By Lemma 2.1 and
the discussion following Theorem 1.8, for $i$ large we may apply Theorem 1.6 to conclude
that for any $0<\alpha\le 1$, the Ricci flows on $(M_i,g_i)$, $g_i(t)$,
exists for $t\in(0, T]$ ($T=T(n,\alpha)$) such that $|\op{Rm}(t^{-1}g_i(t))|\le \alpha$. By Lemma 1.11, for $x_i,y_i\in M_i$ with $d_{t^{-1}g_i(t)}(x_i,y_i)\le 1$,
$$|d_{t^{-1}g_i(t)}(x_i,y_i)-d_{t^{-1}g_i}(x_i,y_i)|\leq \Psi(\alpha|n).$$

Because $(M_i,T^{-1}g_i)@>GH>>(N,T^{-1}g_N)$, let $\phi_i: (M_i,T^{-1}g_i)\to (N,T^{-1}g_N)$ be an $\epsilon_i'$-GHA, $\epsilon'_i\to 0$.
From the above, it is clear that $\phi_i: (M_i,T^{-1}g_i(T))\to (N,T^{-1}g_N)$ is an $(1,\epsilon_i'+\Psi(\alpha|n))$-GHA. By now we can apply Theorem 2.2 to get a nilpotent fibration, and the property that $f_i$ is an $(1-\Psi(\epsilon_i|n))$-H\"older map can been
seen from (1.2.3) and Lemma 2.4 (see below).

Finally, by choosing $\alpha_i\to 0$ and a standard diagonal argument, we obtained the desired
(B1) and (B2) with respect to $g_i$.
\qed\enddemo

\vskip4mm

\subhead b. Proof of Theorem C
\endsubhead

\vskip4mm

Let $M_i$ be as in Theorem C, and we will specify $\delta(n)$ later. By Theorem 1.6, given any $0<\alpha<1$, there is $\delta(\alpha)>0$, such that if every (small) $\rho$-ball in all $\tilde M_i$
satisfies the isoperimetric condition, then Ricci flow solution $g_i(t)$
on $M_i$ exists for $t\in (0, T(n, \rho, \alpha)]$ such that $g_i(0)=g_i$, and
$|\op{Rm}(g_i(t))|\le \frac{\alpha}{t}$. By Theorem 1.8, when all points in $\tilde M_i$ are $(\delta(\alpha),\rho)$-Reifenberg points with $\delta(\alpha)$ suitably small, then the
isoperimetric condition in Theorem 1.6 is satisfied.

Fixing $t\in (0,T(n,\rho,\alpha)]$, passing to a subsequence
we have $$M_i(t):=(M_i,g_i(t))@>\operatorname{GH}>>X_{t}.$$
Because $|\op{Rm}(g_i(t))|\le \frac\alpha t<\frac 1t$, we can apply Theorem 1.4
to get a singular fiber bundle, $f_i: M_i(t)\to X_{t}$, for large $i$. We claim that when
$\alpha$ is chosen small depending $n$, $X_t$ is homeomorphic to $X$ by a $(1-\Psi(\delta|n))$-bi-H\"older map $\phi_t: X_t\to X$ (see Lemma 2.3 below). Then $\phi_t\circ f_i: M_i\to X$
is the desired singular fibration in Theorem C.

\proclaim{Lemma 2.3} Let $M_i(t)@>\operatorname{GH}>>X_t$ be as in the above. Then
we can choose $\alpha$ small depending on $n$ such that there is a bi-H\"older homeomorphism,
$\phi_t: X_t\to X$, satisfying that for any $x,y\in X$ with $d(x,y)\le 1$,
$$c(n)^{-1}d_t(\phi_t(x),\phi_t(y))^{1+\Psi}\le d(x,y)\le c(n)d_t(\phi_t(x),\phi_t(y))^{1-\Psi},$$  where $d_t$ is the distance on $X_t$, and $\Psi=\Psi(\alpha|n)$.
\endproclaim

Lemma 2.3 is a consequence of the following property on $M$; which is a minor
modification of one in (\cite{BW}):

\proclaim{Lemma 2.4} Let $(M,g)$ and $(M,g(t))$ be as in Theorem 1.6.
Then we can choose $\alpha$ small depending on $n$ such that for any $x,y\in M$ with $d_{0}(x,y)\le 1$, we have $$c(n)^{-1}d_{t}(x,y)^{1+\Psi}\le d_{0}(x,y)\le c(n)d_{t}(x,y)^{1-\Psi},
\tag{2.5}$$ where $\Psi=\Psi(\alpha|n)$.
\endproclaim

For the convenience of readers, we will give a proof in the Appendix following
\cite {BW}.

\demo{Proof of Lemma 2.3}

Fixing $t\in (0, T(n,\rho,\alpha)]$,
by Lemma 2.4, the identity maps $\op{id}_i:(\tilde{M}_i,\tilde d_{i})\to(\tilde{M}_i,\tilde d_{i,t})$ are bi-H\"older continuous with uniform coefficients for all $i$. Since $(M_i,d_{i})=(\tilde{M}_i,\tilde d_{i})/\Gamma_i$ and $(M_i,d_{i,t})=(\tilde{M}_i,\tilde d_{i,t})/\Gamma_i$, it is easy to see that the identity maps $\op{id}_i:(M_i,d_{i})\to(M_i,d_{i,t})$ are also bi-H\"older with uniform coefficients for all $i$. In particular, $\op{id}_i:(M_i,d_{i})\to(M_i,d_{i,t})$ are equi-continuous and uniformly bounded. Now by a standard diagonal procedure and by taking a subsequence if necessary,  we may assume that
$\lim_{i\to \infty} \op{id}_i=\psi_t: X\to X_t$, which is H\"older continuous. Similarly, we have $\lim_{i\to \infty} \op{id}_i^{-1}= \phi_t:X_t\to X$. Clearly, $\phi_t\circ\psi_t=\op{id}_{X}$ and $\psi_t\circ\phi_t=\op{id}_{X_t}$. Therefore, $\phi_t$ is a bi-H\"older bijection, and the proof is complete.
\qed\enddemo

\vskip4mm

\subhead c. Proof of Theorem D
\endsubhead

\vskip4mm

Let $M$ be a compact $n$-manifold in Theorem D. As seen in the proof of
Theorem C, there is a complete
Ricci flow solution $g_i(t)$ on $M$, $t\in (0,T]$ and $T=T(n,\rho)$, satisfying
$$|\op{Rm}(g(t))|\le \frac 1t.$$
In particular, $T^{-1}g(T)$ has bounded sectional curvature by
one. By Lemma 1.11, $B_{g(T)}(x,\sqrt T)\subseteq B_{g(0)}(x,1)$. Recall a
basic property of Ricci flow:
$$\op{vol}_{g(t)}(B_{g(0)}(x,1))\le e^{\Psi(t|n)}
\op{vol}_{g(0)}(B_{g(0)}(x,1)),$$
because scalar curvature, $R(g(t))\ge -n(n-1)$, cf. \cite{Ha}.
Thus we may assume $\epsilon(n,\rho)$
is small so that $\op{vol}_{T^{-1}g(T)}(B_{T^{-1}g(T)}(x,1))$ small, and therefore $\op{injrad}_{T^{-1}g(T)}(x)<\epsilon(n)$, a constant in
Theorem 1.5; by Cheeger's injectivity radius estimate (\cite{Ch}).
By now we can apply Theorem 1.5 to conclude Theorem D.
\qed

\vskip4mm

\head Appendix
\endhead

\vskip4mm

\demo{Proof of Theorem 1.2 via embedding method \rm (cf. \cite{Ro1})}

Let's first recall the embedding method: take any $\frac 14$-net
in $N$, $\{y_k\}_{k=1}^s$, and a smooth monotone cut off function,
$h: \Bbb R\to [0,1]$ such that $h(-\infty,\frac 15)\equiv 1$ and $h([\frac 9{10},\infty))\equiv 0$. Then $\Phi(y)=(h(d_N(y_1,y),..., h(d_N(y_s,y))): N\to \Bbb R^s$
is an embedding. Because $\Phi(N)$ is a compact submanifold of $\Bbb R^s$,
there is a $\eta(n)>0$ such that $\exp: D^\perp_\eta(\Phi(N))\to \Bbb R^s$ is an embedding.
If $d_{GH}(M,N)<\epsilon(n)$, then for each $y_i$, take $x_i\in M$ such that $d_H(x_i,y_i)<\epsilon$,
where $d_H$ is an admissible metric on the disjoin union of $N\coprod M$,
$d_H(M,N)<\epsilon$. Then $\Psi: M\to \Bbb R^s$,
$\Psi(x)=(h(d(x_i,x))\in \exp(D^\perp_\eta(\Phi(N)))$, and thus $f=\Phi^{-1}\circ P\circ \Psi:
M\to N$ will be desired map. In practice, we will replace $d(x_i,\cdot)$ by a $C^1$-function $\phi_i$, the average of $d(x_i,\cdot)$ over an $\epsilon$-metric ball, and thus
$f$ is $C^1$-smooth.

We point it out that by assuming $|\op{sec}_M|\le 1$, one can replace $\phi_i$ by
a smooth approximation $\psi_i$. For fixing $x\in B(x_i,\frac 14)$, we define a function,
$f: [d(x_i,x)-1,d(x,x_i)+1]\to \Bbb R_+$, by
$$f(t)=\frac 1{\op{vol}(B(0,\epsilon))}\int_{B(0,\epsilon)}|d(x_i,\exp_xu)-t|^2\op{dvol}_u,$$
where $B(0,\epsilon)\subset B(0,\frac 12)\, (\subset T_xM)$, equipped with the pullback metric. Because the injectivity radius at $0$ is at least $\frac 12$, $f$ is smooth,
which is clearly strictly convex. Assume that $f$ achieves a unique minimum at $c_x$.
We define $\psi_i(x)=c_x$, which is the desired smooth
approximation for $d(x_i,\cdot)$.

Finally, observe that if the metric on $M$ has bounded second derivative, then it is
easy to see that $f$ satisfies (1.2.4). Observe that if $g$ has bounded second derivative, then (1.2.4) holds; while the existence of a $C^1$-close metric of
high regularity can be obtained using Ricci flows (cf. \cite{Sh1, 2}).
\qed\enddemo

\demo{Proof of Lemma 2.4} (cf. \cite{BW})

Given a path $c$, let $f(t)=\op{leng}_{g(t)}(c)$, where $g(t)$ is
the Ricci flow solution, $\frac{\partial }{\partial t}g(t)=-2\op{ Ric}(g(t))$. Using that for $0<s\le w\le t$, $|\op{ Rm}(g(w))|\leq \frac{\alpha}{w}$, by variation of
$f(t)$ one derives
$$\left(\frac{t}{s}\right)^{-(n-1)\alpha}\le \frac{d_t(x,y)}{d_s(x,y)}\le \left(\frac{t}{s}\right)^{(n-1)\alpha}.\tag{2.4.1}$$
Combining that $d_0(x,y)\le d_s(x,y)+\Psi\sqrt{s}$ (see (1.10)) with
the left hand inequality in (2.4.1),
we obtain
$$d_0(x,y)\le\left(\frac{t}{s}\right)^{(n-1)\alpha}d_t(x,y)+\Psi\sqrt{s}.\tag{2.4.2}$$
Choosing $s$ by the equation, $\left(\frac{t}{s}\right)^{(n-1)\alpha}d_t(x,y)=\sqrt{\frac{s}{t}}$, from (2.4.2) we derive
$$d_0(x,y)\le\left(1+\Psi\sqrt{t}\right)d_t(x,y)^{\frac{1}{1+2(n-1)\alpha}},$$
i.e., the right hand side of (2.5).

To see the left hand side of (2.5), we first claim that for $d_0(x,y)\le \sqrt{T(n,\rho,\alpha)}<1$ and $t_0=d_0(x,y)^2$, $d_{t_0}(x,y)\le c_0(n)d_0(x,y)$ and
$c_0(n)$ is a constant depending on $n$.

Assuming the claim, given any two points $x,y\in M$ with $d_0:=d_0(x,y)\le 1$,
and let $l=[d_0/\sqrt{t}]+1$. We divide a $g(0)$-minimal geodesic from $x$ to $y$ into $l$ pieces by $x_0=x,x_1,...,x_l=y$ such that $d_0(x_i,x_{i+1})=l^{-1}d_0$ for $i=0,1,..,l-1$.
Then
$$\aligned
d_0&=\sum_id_0(x_i,x_{i+1})\\&\ge c_0^{-1}\sum_i d_{ l^{-2}d_0^2}(x_i,x_{i+1}) \qquad \text{(applying the claim to $x_i,x_{i+1}$)}\\&\ge c_0^{-1} d_{l^{-2}d_0^2}(x,y)\qquad \text{(a triangle inequality)}\\
&\ge c_0^{-1}d_t(x,y)\left(\frac{l^{-2}d_0^2}{t}\right)^{(n-1)\alpha} \qquad \text{(by (2.4.1) and $t\ge l^{-2}d_0^2$)}\\
&\ge c_1d_0^{2(n-1)\alpha}d_t(x,y). \qquad \text{(by $l^2t\le4$)}\endaligned$$
For $2(n-1)\alpha<1$, we conclude that $d_0\ge c_2(n)d_t(x,y)^{\frac{1}{1-2(n-1)\alpha}}$.


We now verify the claim. Putting $d_{t_0}(x,y)=Ld_0$, without loss of generality
we may assume that $L>>1$. Let $\gamma$ be a $g(0)$-minimal geodesic from $x$ to $y$, and let us choose a partition points $x_1,x_2,\dots ,x_k$ on $\gamma$, where $k=[\frac{L}{2}]$, such that $d_{t_0}(x_i,x_j)>2d_0,i\neq j$. We shall bound on $k$ and thus on $L$ with
a desired constant.

For any $p\in B_{g(t_0)}(x_i,d_0)$,
$$d_0(x_i,p)\le d_{t_0}(x_i,p)+\Psi(\alpha|n)\sqrt{t_0}\le d_0+\Psi(\alpha|n)\sqrt{t_0}=d_0(1+\Psi(\alpha|n)),$$
we derive $$B_{g(t_0)}(x_i,d_0)\subset B_{g(0)}(x_i,2d_0)\subset B_{g(0)}(x,3d_0),$$ as $\alpha$ is sufficient small. Together with $B_{g(t_0)}(x_i,d_0)\cap B_{g(t_0)}(x_j,d_0)=\emptyset$, $i\ne j$,
we derive
$$\split
c(n)d_0^n & \ge\op{vol}_{g(0)}(B_{g(0)}(x,3d_0)) \qquad\text{(by volume comparison)}\\
& \ge\sum_{i=1}^k\op{vol}_{g(0)}(B_{g(t_0)}(x_i,d_0))\\
& \ge e^{-n(n-1)t_0}\sum_{i=1}^k\op{vol}_{g(t_0)}(B_{g(t_0)}(x_i,d_0)) \qquad \text{(by variation of volume)}\\
 &\ge e^{-n(n-1)t_0}k\cdot C(n)d_0^n, \qquad \text{(by Theorem 1.6) }
\endsplit$$
where the last inequality uses
$\op{vol}_{g(t)}(B_{g(t)}(x,\sqrt{t}))\ge C(n)\sqrt{t}^n$ in Theorem 1.6.
By now we see that $k=[L/2]\le c_3(n)$.
\qed\enddemo

\enddocument